\documentclass[11pt]{article}

 \usepackage{hyperref}
\usepackage{graphicx}
\usepackage{amssymb}
\usepackage{epstopdf}
\usepackage{amsmath}

\newcommand{\ggamma}{\boldsymbol{\gamma}}

\newcommand{\flambda}{\frac{1}{\lambda}}
\newcommand{\ex}{e^{i\theta}}
\newcommand{\mex}{e^{-i\theta}}
\newcommand{\zbar}{\bar{z}}

\newcommand{\fracpar}[2]{\frac{\partial #1}{\partial #2}}
\newcommand{\be}{\begin{equation}}
\newcommand{\ee}{\end{equation}}
\newcommand{\bea}{\begin{eqnarray}}
\newcommand{\eea}{\end{eqnarray}}
\newcommand{\eps}{\epsilon}
\newcommand{\C}{\mathbb{C}}
\newcommand{\R}{\mathbb{R}}

\newtheorem{theorem}{Theorem}[section]
\newtheorem{lemma}[theorem]{Lemma}

\newtheorem{proposition}[theorem]{Proposition}

\newenvironment{proof}[1][Proof]{\begin{trivlist}
\item[\hskip \labelsep {\bfseries #1}]}{\end{trivlist}}
\newenvironment{definition}[1][Definition]{\begin{trivlist}
\item[\hskip \labelsep {\bfseries #1}]}{\end{trivlist}}

\newenvironment{remark}[1][Remark]{\begin{trivlist}
\item[\hskip \labelsep {\bfseries #1}]}{\end{trivlist}}

\newcommand{\qed}{\nobreak \ifvmode \relax \else
      \ifdim\lastskip<1.5em \hskip-\lastskip
      \hskip1.5em plus0em minus0.5em \fi \nobreak
      \vrule height0.75em width0.5em depth0.25em\fi}

\title{The Inversion of Ray Transforms on a Conformal Class of Curves}
\author{Nicholas Hoell \footnote{Department of Applied Physics and Applied Mathematics, Columbia University in the City of New York.  Partly supported by NSF research grant number DGE-0221041}  and Guillaume Bal}
\date{May 24, 2010}
\begin{document}
\maketitle
\begin{abstract} \noindent We introduce a technique for recovering a sufficiently smooth function from its ray transform over a wide class of curves in a general region of Euclidean space.  The method is based on  a complexification of the underlying vector fields defining the initial transport and recasting the problem in terms of complex-analytic function theory.  Explicit inversion formulae are then given in a unified form.  The method is then used to give inversion formulae for the attenuated ray transform. 
\end{abstract}

\section{Introduction}
In several engineering situations one deals with data consisting of the line integral of a function and the goal is often to recover that source function from its integral over a class of lines.  In the arena of medical imaging, this arises in positron emission tomography (PET), single photon emission tomography (SPECT), and (originally) CT-scan tomography \cite{natterer2}.  In other applications (in geophysics \cite{shara} and non-destructive electrical imaging techniques \cite{leon, berenstein} such as electrical impedence tomography, EIT) the line integral is instead taken over a class of one-dimensional curves in either Euclidean space or more generally, a Riemannian manifold.  This type of data is referred to as a ray transform.  In geophysics for instance, the problem can arise as the linearization of determining geophysical properties of the Earth from travel-time measurements \cite{shara}.  Quite often the physics will also dictate that the signal suffers some absorption along its trajectory and is attenuated, the data then called, not surprisingly, the attenuated ray transform.  
 
The mathematical applications, properties, and uses of these integral transforms and their inverses are discussed in great detail in  \cite{leon,helgason,helgason1,shara} and include harmonic analysis, algebraic curves, tensor geometry, and partial differential equations to name a few.  Generally, \textit{explicit} inversion formulae over curves other than lines (geodesics of a Riemannian manifold \cite{lee}, say) tend to restrict focus to manifolds with a strong amount of symmetry\cite{helgason, helgason1, helgason2, berenstein, rubin} and do not include the effects of absorption encountered during propagation.  An exception to this statement can be found in \cite{Uhlmann, Uhlmann3}.   We restrict our attention in this paper to curves in a $2$-dimensional region of space.

The method used in this paper generalizes a technique first used in \cite{novikov} for lines in Euclidean space and later generalized in \cite{Bal1} for geodesic rays in hyperbolic geometry for giving an explicit inversion formula for the attenuated ray transform in each case.   The technique we present rests on the complexification of a certain class of differential operators in $\R^2$ which allows us to recast the problem in terms of complex analysis in the unit disc.  Once the problem is cast in this light, we use the classical Poisson formula \cite{fokas} to give us a reconstruction formula.  Excellent introductions to complex analysis and conformal mappings are \cite{sarason, krantz, fokas} and the classic \cite{ahlfors1}.  Good introductions to quasiconformal mappings and Beltrami equations (and their generalizations) can be found in \cite{ahlfors,renelt}. References on Blaschke products and multivalent mappings can be found in \cite{blaschke, garnett}.

An outline of the paper is as follows.  In section \ref{prelim}, the general setup, notation, and a quick review of the essential operators used throughout the paper are presented, together with the main result we are seeking to establish.  In section \ref{complex} we begin the complexification procedure by introducing a new (complex) parameter $\lambda$ into the transport equation introduced in section \ref{prelim} and give a classification of the vector fields under consideration as those of \textbf{type H}.  Much of the heavy lifting is done in the more technical section \ref{solving} where we find and analyze the Green's function of the new parameterized complex partial differential transport equation.  We will establish that \textbf{condition H} is sufficient to guarantee holomorphicity of the solution of this equation in terms of the new parameter $\lambda$.  We evaluate the asymptotics of the solution as our complex parameter $\lambda$ tends to the unit circle from both inside and outside, i.e. as $|\lambda|\to 1^\mp$ and see that in fact its imaginary part depends on the data we are interested in.  Once this is established, we use this fact in section \ref{inv} to give our desired reconstruction formula in the non-attenuated case.  The rest of section \ref{inv} uses the non-attenuated formula to give an integrating factor solution for the attenuated case, which requires an additional constraint to \textbf{condition H}. We offer brief concluding remarks in section \ref{conc}.

\section{Preliminaries}
\label{prelim}
\subsection*{The Stationary Transport Equation}
We let $\ggamma: \R^2 \ni (t,s) \mapsto \ggamma(t,s) \in \Omega \subset \R^2$ be a real-analytic diffeomorphism where $\Omega$ is an open, bounded, simply-connected region of the plane (a domain).  We consider $\R^2 \cong \C$ by the standard isomorphism so that $\ggamma$ is identified with $\gamma^1(t,s)+i\gamma^2(t,s)$. Then, $(w,\bar{w})$ are (independent) coordinates on $\Omega$ where $w\doteq \ggamma(t,s)$. Because $\ggamma$  is a diffeomorphism, its differential is injective and therefore induces a vector field on $\Omega$ via its differential under the rule $(\phi_* X)(f)=X(\phi^* f)$.  Consider $\ggamma_* \fracpar{}{t}$. This gives a non-degenerate field of the following type;
\be 
X|_w=\mu(w)\fracpar{}{w} +\bar{\mu}(w)\fracpar{}{\bar{w}} \qquad w \in \Omega, \qquad |\mu|>0 \nonumber  
\ee 
which acts on pushforwards in $w$ of functions on $\Omega$ and the non-degeneracy is ensured by the regularity of the curves $\ggamma(t,s)$.  The equation of interest is the stationary transport boundary value problem $X|_w u(w)=f(w)$, for $w \in \Omega$, $f(w) \in C^\infty _0(\Omega)$ with $\lim_{t\searrow -\infty} u(w(t,s))=0$, i.e. the BVP
\begin{align}
\mu(w)\fracpar{u}{w}+\bar{\mu}(w)\fracpar{u}{\bar{w}}&=f(w), \qquad w \in \Omega\\
\left.u  \right|_{\partial_- \Omega}&=0
\end{align}
Anticipating our desire to complexify the above equation, we be wanting to exploit the $SO(2)$ symmetry of the unit disc which is a priori not available to us in this more general domain  by appealing to the Riemann mapping theorem \cite{krantz}.  Denote the unit disc by $D^+ \doteq \{z \in \C: |z|<1\}$, the unit circle by $T\doteq \{z \in \C: |z|=1\}$, and $D^-\doteq \C/\{D^+ \cup T\}$. Let $K(z,w)$ be the Bergman function of the domain $\Omega$, then the function
\be
z(w)=\int^w_\zeta \sqrt{\frac{\pi}{K(\zeta,\zeta)}}K(t,\zeta) dt \qquad \zeta \in \Omega
\ee
gives the unique biholomorphism mapping $\Omega$ into $D^+$, the unit disc, with $z(\zeta)=0$,  $z'(\zeta)>0$ as in \cite{nehari} and (t,s)  give coordinates on $D^+$ through composition since $\ggamma^* z$ maps $\R^2$ into $D^+$. Because of this equivalence between our initial domain $\Omega$ and the unit disc all further results will be presented in the disc.  If $\Omega$ was all of $\R^2$ (and the Riemann map was consequently unavailable) the method below will still work since $\R^2$ has the needed rotational symmetry.  

We therefore use $(z,\bar{z})$ as coordinates on $D^+$ and have a new vector field on $D^+$ given by $X|_z =z_*X|_{z(w)}$. and $\mu \to \{z_*\mu\} \fracpar{z}{w} \circ z^{-1}$ and likewise for $\bar{\mu}$. By slight abuse of notation we denote $\{z_*\mu\} \fracpar{z}{w} \circ z^{-1}$ by $\mu(z)$ and $\{z_*\bar{\mu}\} \fracpar{\bar{z}}{\bar{w}} \circ z^{-1}$ by $\bar{\mu(z)}$ so that field of interest is 

\[X|_z=\mu(z)\fracpar{}{z} +\bar{\mu}(z) \fracpar{}{\bar{z}}  ,\qquad z\in D^+\]
We define $t(z) =z_*w_*t$ and $s(z) =z_*w_*s$, smooth functions on $D^+$.  

The method of characteristics gives the following solution to the transport equation $X|_zu(z)=f(z)$, $u(z(-\infty,s))=0$  as 
\be
u(z)=(D_1 f)(z)\doteq \frac{1}{2}\int_\R f(z(t_0,s))sign(t(z)-t_0)dt_0
\ee
and since the integral curves of $X|_z$ are just the image of integral curves, i.e. $\ggamma^* z^*=(z \circ \ggamma)^*$. we define the ray transform of a source function $f(z)$ over the integral curves of $X|_z$ indexed by $s$ to be
\be
(If)(s)=\int_\R f(z(t,s))dt
\ee

We will later be using the following extensions of these operators given below: \newline
\textbf{Symmetrized beam transform}
\[(D_\theta \psi )(z)\doteq \frac{1}{2}\int_\R \psi (\ex z(t_0,s(z \mex )))sign(t(z\mex)-t_0)dt_0 \qquad \psi \in L^1(D^+) \]
\textbf{Ray transform}
\[(I\psi)(s,\ex)=(I_\theta \psi)(s)\doteq \int_\R \psi (\ex z(t,s))dt \qquad \psi \in L^1(D^+)\]
We will always use $\theta$ and $\ex$ interchangeably, its meaning clear from context.

We will have occasion to use the \textbf{Hilbert transform} $H$ of a function defined (see e.g. \cite{stein2}) as the following Calder\'on-Zygmund principal value integral operator 
\be
(H \psi)(x)=\frac{1}{\pi} p.v. \int_\R \frac{\psi(y)}{x-y}dy \qquad \psi \in L^p(\R), \qquad p>1
\ee

Lastly, we will be using the standard \textbf{Poisson kernel} of the unit disc given by $P(z,\theta)=\frac{1-|z|^2}{|1-\mex z|^2}$.  We recall (see, e.g. \cite{evans, taylor}) that the Poisson kernel generates harmonic solutions $v(z)$ of the BVP 
\begin{align}
\Delta v&=0 \qquad z\in D^+ \nonumber \\
\left. v  \right|_{T}&=g \nonumber
\end{align}
given by $v(z)\doteq \frac{1}{2\pi}\int_T P(z,\theta)g(\ex)d\theta$.  The Poisson kernel is also deeply connected to the study of inner functions c.f. \cite{garnett, kumar}.
 \newline
\linebreak
\noindent The main purpose of this article will be to show that given suitable conditions on $\mu(z,\zbar)$ and $s(z,\zbar)$ that 
\[f(z)=\frac{1}{4\pi} \int_0^{2\pi}P(\lambda_i,\theta)X^\bot_\theta H(I_\theta f)(s(z\mex),\ex))d\theta \qquad i=1,...,n \]
where the $\lambda_i(z)$ are functions to be introduced later.  With this established, the above formula is used to give an integrating factor method to find a similar reconstruction formula for the attenuated ray transform along the same lines.  The above is a type of inversion formula known as a filtered backprojection type \cite{natterer2}. 
The procedure used to derive the above main result can be best thought of in the following heuristic scheme
\begin{enumerate}
\item \textbf{Model:}  Writing down the linear stationary transport equation for the dynamics
\item \textbf{Symmetrizing:} Introducing a rotation parameter $\lambda=\ex$ into the integral curves of the transport PDE
\item \textbf{Symmetry-Breaking:} Complexifiying the parameter introduced in step 2 by moving $\lambda$ ``off-shell", i.e. $|\lambda| \neq 1$
\item \textbf{Analysis and Asymptotics:} Evaluating the dependence of solutions to the complexified equation on our parameter $\lambda$ and examining limiting behavior
\item \textbf{Reconstruction:} Using holomorphicity of the solutions to write the inversion formulae as Poisson integrals of their asymptotic boundary values found in step 4
\end{enumerate}
The reader may find some benefit from keeping the above rough outline in mind throughout the following.  In this section, we have finished step 1.  Steps 2 and 3 are handled in the next section.  Step 4 is done in the more technical section \ref{solving}, and the final step is given in section \ref{inv}.

\section{Complexification of the Transport Equation}
\label{complex}
Since $D^+$ is acted on transitively by $SO(2)$  we will define the conformal map $\lambda:(z,\bar{z}) \to (\lambda z, \flambda \bar{z})$, for $\lambda \in T$ the unit circle.  Notice that if $\Phi(\cdot,s)$ is a set of integral curves of $D^+$, that $z^{-1}(\lambda^* \Phi(\cdot,s))$ are conformally related curves in $\Omega$.
Then, for $\lambda \in \{D^+ \cup D^-\}/\{0,\infty\}$ we consider $\lambda_*X|_z \doteq X_\lambda$ to be the so-called ``complexification" of $X|_z$. We remark that $\lambda_*X|_z$ takes the form $\mu(\frac{z}{\lambda},\lambda \zbar)\lambda \fracpar{}{z}+\bar{\mu}(\frac{z}{\lambda},\lambda \bar{z}) \flambda \fracpar{}{\bar{z}}$ or

\be
\label{transport}
X_\lambda =  \xi(z,\lambda) \fracpar{}{z} +\rho(z,\lambda) \fracpar{}{\bar{z}}  \qquad \lambda \in D^\pm /\{0,\infty\} 
\ee
with $\flambda \xi(z,\lambda)= \mu(z,\lambda) \doteq \lambda_*\mu(z)$ and $\lambda \rho(z,\lambda)=\bar{\mu}(z,\lambda)=\lambda_*\bar{\mu}(z)$.  We also define $X_\lambda^{\bot}=\pm i( -\xi(z,\lambda) \fracpar{}{z} +\rho(z,\lambda) \fracpar{}{\bar{z}})$ as a vector field orthogonal to $X_\lambda$ when $\lambda =\ex$.  Namely, $X_\theta \cdot X^\bot_\theta=\pm(\xi(z,\ex),\rho(z,\ex))\cdot (-i\xi(z,\ex),i\rho(z,\ex))=\pm i(|\xi(z,\ex)|^2-|\rho(z,\ex)|^2)=0$ in the standard inner product $(\cdot,\cdot):\C^2 \to \C$. The factor of $i$ is needed to make $X^\bot_\theta u$ real-valued and the choice of $\pm$ is determined by whichever satisfies the condition $X_1^\bot s>0$. Since $X_1^\bot=a(z)z_*\fracpar{}{s}$ for some real-valued $a(z)$, this determines $X^\bot_1$ uniquely.  Since we could just as well reparameterize with $-s$ we will, without any loss of generality, avoid keeping track of signs by just assuming that $X_\lambda^\bot=i(-\xi(z,\lambda)\fracpar{}{z}+\rho(z,\lambda)\fracpar{}{\zbar})$.  

We likewise define $s(z,\lambda)$ and $t(z,\lambda)$ as $\lambda_*s(z)$ and $\lambda_*t(z)$ respectively for $\lambda \in D^\pm/\{0,\infty \}$.  A word on notation: $\fracpar{k}{z}$ and $k_z$ are equivalent, as are $\fracpar{k}{\zbar}$ and $k_{\zbar}$, and we will use them interchangably.  

We remark that equation (\ref{transport}) \textit{has no direct physical meaning} since the complex parameter $\lambda$, when taken to lie away from $T=\partial D^+$,  is in some sense artificial and may be best thought of as a complex parameter indexing a class of complex partial differential equations given in (\ref{transport}).  

Next we reduce the scope of our consideration to the class of vector fields $X_\lambda$ so constructed to consist only of those of \textbf{type H}:
\begin{definition}{A complexified vector field $X_\lambda=a(z,\lambda)\fracpar{}{z}+b(z,\lambda)\fracpar{}{\zbar}$, induced in the manner above as $\lambda_*X|_z$, $\lambda\in D^\pm/\{0,\infty\}$ from a real field $X|_z$, is said to be of \textbf{type H} if the following holds: \begin{itemize}
\item $a(z,\cdot)$ is a holomorphic function of $\lambda$ for $\lambda \in D^+$ and has at least one zero  $\lambda=\lambda_i(z)\in D^+$
\item $b(z,\cdot)$ is a meromorphic function of $\lambda$ for $\lambda \in D^+$ and has no zeroes in $D^+$
\item $\frac{a(z,\cdot)}{b(z,\cdot)}$ is a holomorphic function of $\lambda$ for $\lambda \in D^+$ and has at least one zero  $\lambda=\lambda_i(z)\in D^+$
\item $s(z,\cdot)$, $\fracpar{s(z,\cdot)}{z}$, $\fracpar{s(z,\cdot)}{\zbar}$ are meromorphic functions of $\lambda$ for $\lambda \in D^\pm$
\end{itemize}
where, as in the above, $s(z,\lambda)=\lambda_*s(z)$ is the complexified transverse foliation parameter of the integral curves of $X_\lambda$.
}
\end{definition}

We are, in the above, treating $z$ and $\lambda$ as independent variables, and holomorphicity is to be thought of in the standard way of functions of several complex variables \cite{hormander}.  We are therefore \textit{not} requiring any of the above functions to be holomorphic in the $z$ variable.  

Because $\frac{\xi(z,\lambda)}{\rho(z,\lambda)}$ is holomorphic in $\lambda \in D^+$  its zeroes are isolated. Also, since $\frac{\xi(z,\lambda)}{\rho(z,\lambda)}$ is holomorphic for $\lambda \in D^+$ and since conformal mappings map boundaries of Jordan domains into boundaries of Jordan domains (see \cite{nehari}), then $\frac{\mu(z,\ex)}{\bar{\mu}(z,\ex)}=\frac{\mu(y)}{\bar{\mu}(y)}$ for some $y \in T$ and thus $ |\left. \frac{\xi(z,\lambda)}{\rho(z,\lambda)} | \right |_{|\lambda|=1}=1$.  Since we assumed that there is at least one zero $\lambda_i$, the maximum principle implies that $|\frac{\xi(z,\lambda)}{\rho(z,\lambda)} |<1$ for $\lambda \in D^+$.  
We then get the following simple 
\begin{lemma}{$\frac{\xi(z,\lambda)}{\rho(z,\lambda)}$  has a finite number of zeros, $\lambda_i(z)$ with multiplicities $m_i(z)$}
\begin{proof}
This is a simple consequence of the argument principle (\cite{krantz}).  Namely, one has
\[\sum_i m_i =\frac{1}{2\pi i} \int_{|\lambda|=1} \frac{\fracpar{}{\lambda}{\frac{\xi(z,\lambda)}{\rho(z,\lambda)}} }{\xi(z,\lambda)}\rho(z,\lambda)d\lambda \]
and $\frac{\xi(z,\lambda)}{\rho(z,\lambda)}$ is holomorphic, hence so is $\fracpar{}{\lambda}{\frac{\xi(z,\lambda)}{\rho(z,\lambda)}}$, on the region $D^+$. They are also both continuous on $T$.  Therefore, $|\fracpar{}{\lambda}{\frac{\xi(z,\lambda)}{\rho(z,\lambda)}}|<M<\infty$ for $\lambda \in \overline{D^+}$. Thus, 
\be
\sum_i m_i \leq \frac{1}{2\pi} |\int_{|\lambda|=1} \frac{\fracpar{}{\lambda}{\frac{\xi(z,\lambda)}{\rho(z,\lambda)}}}{\xi(z,\lambda)}\rho(z,\lambda)d\lambda| < \frac{1}{2\pi} \int_0^{2\pi} M d\theta = M
\ee
\qed
\end{proof}
\end{lemma}

\noindent Henceforth $\lambda_i$ will always be used to indicate a value in the unit disc for which $\frac{\xi(z,\lambda)}{\rho(z,\lambda)}$ (and $\xi$) vanishes. The bounded holomorphic functions mapping the unit disc onto itself and having a finite number of zeroes can be uniquely written as a finite Blashke product (c.f. \cite{blaschke, garnett}) so that $\frac{\xi(z,\lambda)}{\rho(z,\lambda)}$ can be given in the form $\zeta (z)\Pi_{i=1}^n(\frac{\lambda-\lambda_i}{1-\lambda \bar{\lambda}_i})^{m_i}$ with $|\zeta(z)|=1$, and with $m_i$ and $\lambda_i$ possibly depending on $z$.  

Furthermore, since  $|\frac{\xi(z,\lambda)}{\rho(z,\lambda)} |<1$ for $\lambda \in D^+$ we also have that the complexified transport equation $X_\lambda u(z,\lambda)=f(z)$ can be rewritten as
\be
\label{belt}
u_{\zbar}(z,\lambda)=\frac{\xi(z,\lambda)}{\rho(z,\lambda)}u_z (z,\lambda)+\frac{f(z)}{\rho(z,\lambda)} \ee
which is a forced Beltrami equation as in \cite{ahlfors2,begehr}.

\section{Solving the Complexified Equation}
\label{solving}
In trying to solve the complexified transport equation
\be
\label{stationary}
X_\lambda u(z,\lambda)=f(z)
\ee
we will again be changing variables.  Notice that $X_\lambda s(z,\lambda)=0$ on that region so that $s$ is still a constant of the dynamics. This is obvious from the fact that integral curves are mapped by diffeomorphisms to integral curves, however to be precise, when $|\lambda|\neq 0$,  
\bea
X_\lambda s(z,\lambda)&=&\lambda_*X|_z\lambda_*s(z)=\lambda_*z_*w_*\fracpar{}{t}\lambda_*z_*w_*s=(\lambda\circ z\circ w)_*\fracpar{}{t}(\lambda \circ z \circ w)_*s \nonumber \\
&=&(\lambda \circ z \circ w)_* \fracpar{s}{t}=0 \nonumber
\eea
since $s$ and $t$ are independent coordinates. Thus
\be
\xi(z,\lambda)\fracpar{s(z,\lambda)}{z}+\rho(z,\lambda) \fracpar{s(z,\lambda)}{\bar{z}}=0
\ee
The Riemann removable singularities theorem \cite{krantz} applies when $\lambda=0$. We will need the following 
\begin{lemma}{On $0<|\lambda|<1$ the Jacobian $\partial s(z) \doteq |s_z(z,\lambda)|^2-|s_{\zbar} (z,\lambda)|^2$ is positive}
\begin{proof}

\noindent Since $(t,s) \mapsto (z,\zbar)$ is a diffeomorphism and $\lambda: (z,\zbar) \mapsto (\frac{z}{\lambda},\zbar \lambda)$ is conformal on $0<|\lambda|<1$ 
\be
\begin{vmatrix}
\fracpar{s(z,\lambda)}{(\frac{z}{\lambda})} & \fracpar{t(z,\lambda)}{(\frac{z}{\lambda})} \\
\fracpar{s(z,\lambda)}{(\zbar \lambda)} & \fracpar{t(z,\lambda)}{(\zbar \lambda)} 
\end{vmatrix}
\neq 0
\ee
so that 
\[|\fracpar{s(z,\lambda)}{z}\fracpar{t(z,\lambda)}{\zbar}-\fracpar{s(z,\lambda)}{\zbar}\fracpar{t(z,\lambda)}{z}| \leq |s_z(z,\lambda)|( |t_{\zbar}(z,\lambda)|+|\frac{\xi(z,\lambda)}{\rho(z,\lambda)} || t_z (z,\lambda)|) \]
implies $ |s_z(z,\lambda)|^2 \neq 0$.  Then, \[\partial s(z)= |s_z(z,\lambda)|^2-|\frac{\xi(z,\lambda)}{\rho(z,\lambda)} s_z(z,\lambda)|^2 \geq  |s_z(z,\lambda)|^2(1-|\frac{\xi(z,\lambda)}{\rho(z,\lambda)}|^2)>0\]
since  $|\frac{\xi(z,\lambda)}{\rho(z,\lambda)} |<1$ for $\lambda \in D^+$.
\qed
\end{proof}
\end{lemma}

Since $X_\lambda s(z,\lambda)=0$,  $s_*X_\lambda=s_*X_\lambda \bar{s}(z,\lambda) \fracpar{}{\bar{s}}$. We are interested in solving $X_\lambda G_\lambda(z;z_0)=\delta(z-z_0)$ and we can achieve this by solving $s_*X_\lambda \bar{s}(z,\lambda) \fracpar{}{\bar{s}} (s_* G_\lambda) =|\partial s(z)| \delta(s(z,\lambda)-s_0)$. We therefore need to compute the term $s_*X_\lambda \bar{s}(z,\lambda)$. To this end, 
\[ \xi(z,\lambda)\fracpar{s(z,\lambda)}{z}+\rho(z,\lambda)\fracpar{s(z,\lambda)}{\bar{z}}=0\]
implies  \[\xi(z,\lambda)=-\rho(z,\lambda)\frac{ \fracpar{s(z,\lambda)}{\bar{z}}}{\fracpar{s(z,\lambda)}{z}} \]
whence 
\begin{align} 
 \xi(z,\lambda)\fracpar{\bar{s}(z,\lambda)}{z}+\rho(z,\lambda)\fracpar{\bar{s}(z,\lambda)}{\bar{z}} &=
-\rho(z,\lambda)\frac{ \fracpar{s(z,\lambda)}{\bar{z}}}{\fracpar{s(z,\lambda)}{z}} \fracpar{\bar{s}(z,\lambda)}{z} +\rho(z,\lambda)\fracpar{\bar{s}(z,\lambda)}{\bar{z}} \nonumber \\
 &=\frac{\rho(z,\lambda)}{\fracpar{s(z,\lambda)}{z}}(|s_z(z,\lambda)|^2-|s_{\bar{z}}(z,\lambda)|^2) \nonumber \\
 &= \frac{1}{Q(z,\lambda)} \partial s(z) \nonumber 
\end{align}
with $Q(z,\lambda)\doteq \frac{\fracpar{s(z,\lambda)}{z}}{\rho(z,\lambda)}$. By recalling that $|\frac{\xi}{\rho}|>1$ for $|\lambda|>1$ and going through the preceding lemma \textit{mutatis mutandis} we see that $\partial s(z)$ is likewise negative on $D^-$ and hence the Jacobian of $s(z,\lambda)$ \textit{switches sign} when $\lambda \in D^\pm$ so that our fundamental equation then becomes 
\[s_*\frac{1}{Q(z,\lambda)}\fracpar{}{\bar{s}}s_*G_\lambda=sign(1-|\lambda|)\delta(s(z,\lambda)-s(z_0,\lambda))\]
which gives $G_\lambda(z;z_0)=\frac{sign(1-|\lambda|)Q(z_0,\lambda)}{\pi(s(z)-s(z_0))}$ or rather 
\be
\label{greens}
G_\lambda(z;z_0)=\frac{sign(1-|\lambda|)\frac{1}{\rho(z_0,\lambda)}\fracpar{s(z,\lambda)}{z}|_{z_0} }{\pi(s(z,\lambda)-s(z_0,\lambda))} , \qquad \lambda \in D^\pm /\{0,\infty\}
\ee
so that $u(z,\lambda)=\int_{D^+} G_\lambda(z;z_0)f(z_0)d\mu(z_0)$ solves $X_\lambda u(z,\lambda)=f(z)$ for $\lambda \in D^\pm/\{0,\infty\}$.  We have used the fact that $\fracpar{}{z}\frac{1}{\pi \zbar}=\delta(z)$ as shown in \cite{krantz}.  
\newline

\begin{remark}{We will only make use of results in our formula which follow from \textbf{condition H} and thus results like (\ref{greens}) are only used when $\lambda \in D^+$. We will however present many results for $\lambda \in D^-$ with the understanding that given an appropriate generalization of \textbf{condition H} (involving constraint on $\xi$ and $\rho$ for $\lambda \in \C/\bar{D}^+$) the results are true. The advantage to this approach is it makes use of the symmetries and parallels of several of the formulae for $\lambda \in D^\pm$. Thus, in the ``$-$" versions of several results, \textbf{condition H} is necessary but not sufficient.} 
\end{remark}

With the above in mind, since 
 
 \begin{align}
 \frac{\partial(t,s)}{\partial(z,\zbar)}\frac{\fracpar{\zbar}{t}}{\fracpar{s}{z}}&=\fracpar{\zbar}{t}\fracpar{s}{\zbar}\fracpar{t}{z}\frac{1}{\fracpar{s}{z}}-\fracpar{\zbar}{t}\fracpar{t}{\zbar} \nonumber \\
 &=-(\fracpar{z}{t}\fracpar{t}{z}+\fracpar{\zbar}{t}\fracpar{t}{\zbar}) \nonumber \\
&=-z_*\fracpar{t}{t}
 \end{align}
 we can rewrite $G_\lambda(z;z_0)$ as 
 \be
 \label{green}
 G_\lambda(z,z_0)=-\lambda \frac{ \left.  \frac{\partial(t,s)}{\partial(z,\zbar)} \right |_{z_0} }{\pi (s(z)-s(z_0))}
 \ee
 Then for $\psi \in C^\infty_0(D^+)$ and $d\mu (z) =\frac{dz d\zbar}{2i}=dxdy$, the standard Lebesgue measure on $\R^2 \cong \C$ we have 
 \begin{align}
 \int_{\lambda(D^+)} \psi(z_0)\lambda \frac{ \left.  \frac{\partial(t,s)}{\partial(z,\zbar)} \right |_{z_0} }{\pi (s(z)-s(z_0))}d\mu(z_0)  &=\int_{D^+}(\lambda^* \psi) \frac{ \left.  \frac{\partial(t,s)}{\partial(z,\zbar)} \right |_{z_0} }{\pi (s(z)-s(z_0))}d\mu(z_0) \nonumber \\
 &=\frac{1}{2\pi i} \int_\R \int_\R \frac{\lambda^*\psi(z(t_0,s_0)) dt_0 ds_0}{s-s_0} 
 \end{align}
 so that $\int_{D^+} \psi(z_0) G_\lambda (z;z_0) d\mu(z_0)$ stays bounded since $\lambda$ fixes the unit disc.  A similar argument works for $\lambda \in D^-/\{\infty\}$.  Because of the meromorphy assumptions stated in condition H, we have that when $z \neq z_0$, $G_\lambda (z;z_0)$ is a holomorphic function for $\lambda \in D^\pm/\{0,\infty\}$.  Since $s(z,\lambda)$ and $s_z(z,\lambda)$ have the same order of (possible) pole at zero, $G_\lambda(z;z_0)$ stays bounded even at $\lambda =0$ and we get the following
 
 \begin{proposition}{$u(z,\lambda)$ is holomorphic for $\lambda \in D^\pm$}
 \end{proposition} 

 A similar argument applied to $\fracpar{}{z}G_\lambda(z;z_0)$ and $\fracpar{}{\zbar}G_\lambda(z;z_0)$ shows that $u_z (z,\lambda)$ and $u_{\zbar}(z,\lambda)$ respectively are also complex-analytic for $\lambda \in D^\pm$, a fact we will make use of in our final reconstruction formulae.

\subsection*{Boundary Behavior}

We will be using the boundary values $u(z,\lambda)|_{\lambda \in T}$ to arrive at a reconstruction formula.  Therefore, ignoring the signum for the moment and letting  $\psi \in C_0^\infty (D^+)$ be a test function, then, using (\ref{green}), the two-form $\psi(z_0)G_\lambda(z;z_0)d\mu(z_0)$ equals \[-\lambda_*\{\psi(\lambda^*z(t_0,s_0))\frac{1}{2\pi i(s-s_0)} dt_0 ds_0\}\]
so that we get the following

\begin{proposition}{$u_{\pm}(z,\ex) \doteq \lim_{D^\pm \ni \lambda \to \ex}u(z,\lambda) = \mp \frac{1}{2i}(HI_\theta f)(s(e^{-i\theta}z),\theta)+(D_\theta f)(z)$ where the Hilbert transform $H$ is taken with respect to the first variable.} 
\begin{proof}
First we examine $\frac{1}{s(z,\lambda)-s(z_0,\lambda)}$ when $\lambda =1-\eps$ ($\eps <<1$) and use the fact that $s(z,1-\eps)=s(z,1)-\eps s(z,1)+o(\eps^2)$ together with $X_\lambda s(z,\lambda)=0$ to get  
\begin{align}
O(1)&:& (\xi(z,1)\fracpar{}{z}+\rho(z,1)\fracpar{}{\zbar})s(z,1)&=0 \nonumber \\
O(\eps) &:& (\xi(z,1)\fracpar{}{z}+\rho(z,1) \fracpar{}{\zbar})s'(z,1)&=-(\xi '(z,1)\fracpar{}{z} +\rho '(z,1)\fracpar{}{\zbar})s(z,1) \nonumber
\end{align}
and 
\begin{align}
-(\xi '(z,1)\fracpar{}{z} +\rho '(z,1)\fracpar{}{\zbar})s(z,1) &= -(\xi'(z,1)-\rho'(z,1)\frac{\xi(z,1)}{\rho(z,1)})s_z(z,1) \nonumber \\
&=-\xi(z,1)s_z(z,1)\{ \frac{\xi'(z,1)}{\xi(z,1)}-\frac{\rho'(z,1)}{\rho(z,1)}\} \nonumber \\
&=-\xi(z,1)s_z(z,1)\frac{  (\left. \fracpar{}{\lambda} \frac{\xi}{\rho}) \right |_{\lambda =1} }{\frac{\xi(z,1)}{\rho(z,1)}} 
\end{align}
so that
\[X_1 i s'(z,1)=-i\xi(z,1)s_z(z,1)\frac{  (\left. \fracpar{}{\lambda} \frac{\xi}{\rho}) \right |_{\lambda =1} }{\frac{\xi(z,1)}{\rho(z,1)}} \]
By a similar argument one can show
\[X_1 i s'(z,1)=i \rho(z,1) s_{\zbar}(z,1)\frac{  (\left. \fracpar{}{\lambda} \frac{\xi}{\rho}) \right |_{\lambda =1} }{\frac{\xi(z,1)}{\rho(z,1)}} \]
so that
\be
\label{mess}
X_1 i s'(z,1)=\frac{1}{2} \frac{  (\left. \fracpar{}{\lambda} \frac{\xi}{\rho}) \right |_{\lambda =1} }{\frac{\xi(z,1)}{\rho(z,1)}}X^\bot_1 s(z,1) 
\ee
Since $\frac{\xi}{\rho}$ is given as a finite Blashke product $\zeta(z)\Pi_{i=1}^n(\frac{\lambda-\lambda_i(z)}{1-\lambda \bar{\lambda}_i})^{m_i(z)}$, we see that $\frac{\fracpar{}{\lambda}\frac{\xi(z,\lambda)}{\rho(z,\lambda)}}{\frac{\xi(z,\lambda)}{\rho(z,\lambda)}}=\sum_{j > 0 }m_j \frac{1-|\lambda_j|^2}{(\lambda-\lambda_j)(1-\bar{\lambda}_j \lambda)}$ so that $ \frac{  (\left. \fracpar{}{\lambda} \frac{\xi}{\rho}) \right |_{\lambda =1} }{\frac{\xi(z,1)}{\rho(z,1)}} >0$, which, when combined with $X^\bot_1s(z,1)>0$ gives from (\ref{mess}) that  
\[X_1 i s'(z,1) >0\]
and therefore \[sign(is'(z,1)-is'(z_0,1))=sign(t(z,1)-t(z_0,1)).\]  Then we look at 
\be 
\int_{D^+} G_{1-\eps}(z;z_0)\psi(z_0)d\mu(z_0) \to -\frac{1}{2\pi i} \int_\R \int_\R \frac{\psi(z(t_0,s_0))}{s(z,1-\eps)-s(z_0,1-\eps)}ds_0dt_0\nonumber
\ee
 which as we have shown is 
 \bea 
&& -\frac{1}{2\pi i} \int_\R \int_\R \frac{\psi(z(t_0,s_0))}{s(z,1)-s(z_0,1)-\eps(s'(z,1)-s'(z_0,1))}ds_0dt_0 \to \nonumber \\ 
 && -\frac{1}{2\pi i} \int_\R \int_\R \frac{\psi(z(t_0,s_0))}{s(z,1)-s(z_0,1)}ds_0dt_0 \nonumber \\
 && +\frac{1}{2} \int_\R \int_\R \delta(s(z,1)-s(z_0,1)) sign(is'(z,1)-is'(z_0,1))\psi(z(t_0,s_0))dt_0ds_0 \nonumber \\
 &=&\frac{-1}{2i}H(I_\theta \psi)(s(z),1) + \frac{1}{2}  \int_\R sign(t(z,1)-t_0)\psi(z(t_0,s(z,1)))dt_0 
  \eea
  Thus, we have 
  \be
u_+(z,1)=\frac{-1}{2i}H(I_\theta \psi)(s(z),1)+(D_1 \psi)(z)
\ee
For the general case, $G_{\ex}(z;z_0)=G_1(\mex z;\mex z_0)$ shows
\be
u_+(z,\ex)=\frac{-1}{2i}H(I_\theta \psi)(s(z\mex),\ex)+(D_\theta \psi)(z).
\ee
An identical argument for $u_-(z,\ex)$ shows that
\be
 u_\pm(z,\ex)=\mp\frac{1}{2i}H(I_\theta \psi)(s(z\mex),\ex)+(D_\theta \psi)(z)
\ee
\qed
\end{proof}
\end{proposition}

\section{Inversion Formulae}
\label{inv}
\subsection{No Attenuation} 
We can now prove our main result.
\begin{theorem}{If $X_\lambda$ is a vector field of \textbf{type H}, $\left. \frac{\xi(z\lambda)}{\rho(z\lambda)}\right |_{\lambda_i(z)}=0$ for $i=1,...,n$ and $f(z) \in C^\infty_0(D^+)$, then 
\[f(z)=\frac{1}{4\pi} \int_0^{2\pi}P(\lambda_i,\theta)X^\bot_\theta H(I_\theta f)(s(z\mex),\ex)d\theta \] gives an exact reconstruction formula for the density $f$ based on the data $I_\theta f$ of ray transforms of $f$ over the integral curves of $X_\theta$.}

\begin{proof}
With  $P(z,\theta)=\frac{1-|z|^2}{|1-\mex z|^2}$, the Poisson kernel of the unit disc, and Cauchy's formula for holomorphic functions. one has (\cite{fokas}) that
\begin{align}
\label{X1}
X_{\lambda_i}u(z,\lambda_i)&=\frac{i}{4\pi} \int_0^{2\pi}P(\lambda_i,\theta)X_\theta H(I_\theta f)(s(z\mex),\ex)d\theta \nonumber \\
 &+ \frac{1}{2\pi} \int_0^{2\pi}P(\lambda_i,\theta)X_\theta (D_\theta f)(z)d\theta
\end{align}
 so that
 
 \[X_{\lambda_i}u(z,\lambda_i)=f(z)+\frac{i}{4\pi} \int_0^{2\pi}P(\lambda_i,\theta)X_\theta H(I_\theta f)(s(z\mex),\ex)d\theta \]
 whereas 

\begin{align}
 \label{X2}
 X^\bot_{\lambda_i}u(z,\lambda_i)&=\frac{i}{4\pi} \int_0^{2\pi}P(\lambda_i,\theta)X^\bot_\theta H(I_\theta f)(s(z\mex),\ex)d\theta \nonumber \\
  &+ \frac{1}{2\pi} \int_0^{2\pi}P(\lambda_i,\theta)X^\bot_\theta (D_\theta f)(z)d\theta
 \end{align}
Then since $X_\lambda=\xi(z,\lambda) \fracpar{}{z} +\rho(z,\lambda) \fracpar{}{\bar{z}}$, $X_\lambda^{\bot}= i( -\xi(z,\lambda) \fracpar{}{z} +\rho(z,\lambda) \fracpar{}{\bar{z}})$ and $\xi(z,\lambda_i)=0$, we have that  \[i X_{\lambda_i}u(z,\lambda_i)=X^\bot_{\lambda_i}u(z,\lambda_i)\] so that, on equating real and imaginary parts of (\ref{X1}) and (\ref{X2}), we get

 \be
 \frac{1}{2\pi} \int_0^{2\pi}P(\lambda_i,\theta)X^\bot_\theta (D_\theta f)(z)d\theta=-\frac{1}{4\pi} \int_0^{2\pi}P(\lambda_i,\theta)X_\theta H(I_\theta f)(s(z\mex),\ex)d\theta \nonumber
 \ee
 and
  \be
 \label{recon}
 f(z)=\frac{1}{4\pi} \int_0^{2\pi}P(\lambda_i,\theta)X^\bot_\theta H(I_\theta f)(s(z\mex),\ex)d\theta
 \ee
 \qed
 \end{proof}
 \end{theorem}
 
It's clear that formula (\ref{recon}) could just as well be written in terms of the jump function (from the viewpoint of $D^\pm$)
\[\phi(z,\ex) \doteq u_+(z,\ex)-u_-(z,\ex)=iH(I_\theta f)(s(z\mex),\ex) \]
as 
\be
f(z)=\frac{1}{4\pi} \int_0^{2\pi}P(\lambda_i,\theta)X^\bot_\theta(-i\phi(z,\ex))d\theta
\ee
an observation which will be useful in the next section.  Recalling our previous \textbf{remark} about using only results from $D^+$ we could just as well use
\[\phi(z,\ex) \doteq 2i\Im(u_+(z,\ex)) \]
and remember that invoking $D^-$ is only a useful mnemonic.  

\subsection{Attenuated Ray Transform and Inversion Formulae}

We add a real-valued attenuation term $a(z)\in C_0^\infty (D^+)$ to the complexified stationary transport equation to get

\be
(X_\lambda + a(z) ) u(z,\lambda) =f(z) \qquad \lambda \in D^\pm
\ee 
Using our Green's function $G_\lambda(z;z_0)$, we define
\be
h(z,\lambda) \doteq \int_{D^+} G_\lambda (z;z_0)a(z_0)d\mu(z_0)
\ee
and we use an integrating factor approach as follows
\[e^{h(z,\lambda)}X_\lambda u(z,\lambda)+e^{h(z,\lambda)}a(z)u(z,\lambda)=e^{h(z,\lambda)}f(z)\]
so that
\[X_\lambda e^{h(z,\lambda)}u(z,\lambda)=e^{h(z,\lambda)}f(z)\]
whence 

\be
u(z,\lambda) =\int_{D^+}  G_\lambda(z;z_0)e^{h(z_0,\lambda)-h(z,\lambda)}f(z_0)d\mu(z_0) 
\ee
Now, since 
\be
h_\pm (z,\ex)= \mp \frac{1}{2i} (HI_\theta a)(s(z \mex), \theta)) + (D_\theta a)(z) \nonumber
\ee
as before, we have the solution of the attenuated transport equation admits the following boundary values as $|\lambda| \to 1^\mp$
\bea
u_\pm (z,\ex)&=&\frac{\mp e^{-h_\pm (z,\ex)}}{2i} [ HI_\theta \{ e^{h_\pm (\cdot,\ex)}f \}(s(z \mex ),\theta) \mp 2i (D_\theta   e^{h_\pm (\cdot,\ex)}f )(z) ] \nonumber \\
&=& \frac{\mp e^{-h_\pm (z,\ex)}}{2i}[ HI_\theta \{ e^{\frac{\mp 1}{2i} (HI_\theta)a(s(\mex \cdot), \theta)}f(\cdot) e^{(D_\theta a)(\cdot)} \}(s(z \mex ),\theta) \nonumber \\
&\mp&  2i (D_\theta  e^{\frac{\mp 1}{2i} (HI_\theta)a(s(\mex \cdot), \theta)}f(\cdot) e^{(D_\theta a)(\cdot)} )(z)] \nonumber
\eea

\noindent Defining
\be
(I_{a,\theta} f)(s) \doteq I_\theta (f(\cdot)  e^{(D_\theta a)(\cdot)} )(s)
\ee
and recalling that $I_\theta$ involves integration in $t$, not $s$ (as does $D_\theta$) and therefore 
\be
u_\pm(z,\ex)=\frac{\mp e^{-h_\pm (z,\ex)}}{2i} H(e^{\frac{\mp 1}{2i} (H (I_\theta a)(s(\mex \cdot), \theta)} I_{a,\theta} f)(s(z \mex),\theta)+e^{-(D_\theta a)(z)} (D_\theta f(\cdot) e^{(D_\theta a)(\cdot)})(z) \nonumber 
\ee

\noindent so that 
\begin{align}
\phi(z,\ex) &\doteq (u_+-u_-)(z,\ex) =-\frac{e^{-h_- (z,\ex)}}{2i} H(e^{\frac{ 1}{2i} H(I_\theta a)(s(\mex \cdot), \theta)} I_{a,\theta} f)(s(z \mex),\theta) \nonumber \\
&  -\frac{ e^{-h_+ (z,\ex)}}{2i} H(e^{-\frac{1}{2i} H(I_\theta a)(s(\mex \cdot), \theta)} I_{a,\theta} f)(s(z \mex),\theta) \nonumber \\
&= -\frac{e^{-(D_\theta a)(z)}}{2i} \{  e^{\frac{1}{2i} H(I_\theta a)(s(z \mex), \theta)} H(e^{\frac{1}{2i} H(I_\theta a)(s(\mex \cdot), \theta)} I_{a,\theta} f) \nonumber \\
&+ e^{-\frac{1}{2i} H(I_\theta a)(s(z \mex), \theta)} H(e^{-\frac{1}{2i} H(I_\theta a)(s(\mex \cdot), \theta)} I_{a,\theta} f) \}(s(z \mex),\theta) \nonumber
\end{align}
We define $C \doteq \cos (\frac{H(I_\theta a)(s(z \mex),\theta)}{2})$ and $S\doteq \sin ( \frac{H(I_\theta a)(s(z \mex),\theta)}{2})$. Then
\bea
\phi(z,\ex)&=& -\frac{e^{-(D_\theta a)(z)}}{2i} [ (C-iS)H \{ (C-iS)I_{a,\theta} f\}  +(C+iS)H \{ (C+iS)I_{a,\theta} f\} ](s(z \mex),\theta) \nonumber \\
&=&ie^{-(D_\theta a)(z)} \Re \{ (C-iS)H [ (C-iS)I_{a,\theta} f] (s(z \mex),\theta) \} \nonumber \\
&=&ie^{-(D_\theta a)(z)}(CH(CI_{a,\theta} f) (s(z \mex),\theta) +SH(SI_{a,\theta} f) (s(z \mex),\theta) ) \\
&\doteq&ie^{-(D_\theta a)(z)}(H_aI_{a,\theta} f) (s(z \mex),\theta) 
\eea
where $H_a :f \mapsto CH(CI_{a,\theta} f) (s(z \mex),\theta) +SH(SI_{a,\theta} f) (s(z \mex),\theta) $. We then can proceed in a manner similar to before since we have that $e^{h(z,\lambda)} u(z,\lambda)$ (along with its derivatives) is holomorphic and solves $X_\lambda e^{h(z,\lambda)} u(z,\lambda)=e^{h(z,\lambda)}f(z)$.  We stipulate, in addition to $X_\lambda$ being of \textbf{type H} that, furthermore, $u(z,\lambda_k)=0$ for all $\lambda_k(z)$ for which $\xi(z,\lambda)=0$.  Under this \textbf{additional assumption}, we see that in fact 
\be
i(X_{\lambda_i} u(z,\lambda_i)+a(z)u(z,\lambda_i)=X^\bot_{\lambda_i}u(z,\lambda_i)
\ee
and we have proven that 
\begin{theorem}{If $X_\lambda$ is a vector field of type H, $u(z,\lambda_i)=0$ and $f \in C^\infty_0(D^+)$, then \[f(z)=\frac{1}{4\pi} \int_0^{2\pi}P(\lambda_i,\theta) X^\bot_\theta (e^{-(D_\theta a)(z)}H_aI_{a,\theta} f) (s(z \mex),\theta)) d\theta \] gives an exact reconstruction formula for the density $f$ based on the data $I_{a,\theta} f$ of attenuated ray transforms of $f$ over the integral curves of $X_\theta$.}
\end{theorem}

\section{Conclusions}
\label{conc}
The method of complexification presented in the preceding allows for a compact unification of the inversion formulae given for ray transforms on both Euclidean space \cite{novikov} and the Poincar\'e hyperbolic disc \cite{Bal1}.  Extending the class of vector fields amenable to the aforementioned scheme beyond those of \textbf{type H}\ remains an open problem.  Since the analyticity properties of the coefficients of the vector fields, ensured by the \textbf{condition H}, were used to justify the holomorphy of the Green's function it is unclear how one could alter the method in the absence of such a condition, although the recent \cite{Uhlmann3} may yield some insight.  With that in mind, there remains the question of finding sufficient (or even necessary) conditions on the \textit{initial} vector field being holomorphic after the complexification used above.  Real-analyticity is perhaps the simplest necessary condition, but presumably there are much more stringent ones.  There also remains the question of when $u(z,\lambda_i)=0$.  

Lastly, we remark that the only symmetry of the equations occurs when $\lambda \in T$, which is \textit{not} where the analysis takes place.  In fact, we must break the symmetry in order to arrive at our solution and find our minima $\lambda_i$.  \textit{Informally}, this procedure is analogous to the so-called Higgs mechanism for gauge-invariant spontaneous symmetry-breaking of a complex scalar field used in the Standard Model of particle physics \cite{griffiths,rubakov}.
\bibliographystyle{siam}

\end{document}